\documentclass[11pt]{article}
\usepackage{graphicx}
\usepackage{amsfonts}
  \def\fH{{\cal H}}

\usepackage{theorem}
\newtheorem{Lem}{Lemma}[section]

\newtheorem{The}[Lem]{Theorem}

\newtheorem{Cor}[Lem]{Corollary}

\newtheorem{Rem}[Lem]{Remark}

\newcommand{\qed}{\hbox{\rule{6pt}{6pt}}}

\begin{document}
\title{Alternative reverse inequalities for Young's inequality}
\author{Shigeru Furuichi$^1$\footnote{E-mail:furuichi@chs.nihon-u.ac.jp} and Nicu\c{s}or Minculete$^2$\footnote{E-mail:minculeten@yahoo.com}\\
$^1${\small Department of Computer Science and System Analysis,}\\
{\small College of Humanities and Sciences, Nihon University,}\\
{\small 3-25-40, Sakurajyousui, Setagaya-ku, Tokyo, 156-8550, Japan}\\
$^2${\small \lq\lq Dimitrie Cantemir\rq\rq  University, Bra\c{s}ov, 500068, Rom{a}nia}}
\date{}
\maketitle
{\bf Abstract.}
Two reverse inequalities for Young's inequality were shown by M. Tominaga, using Specht ratio.
 In this short paper, we show alternative reverse inequalities for Young's inequality without using Specht ratio.
\vspace{3mm}


\vspace{3mm}
\vspace{3mm}

\section{Introduction}
It is well known the Young inequality
\begin{equation} \label{eq_1_1}
a^{1-\lambda}b^{\lambda} \leq   (1-\lambda) a+ \lambda b, 
\end{equation}
for positive real numbers $a,b$ and $\lambda \in [0,1]$. 
See \cite{BK,Dra,KM,Ald2,Ald1,Mit,Furu1,Furu2,Min1,Min2} for improvements of Young's inequality and their recent advances.

One of reverse inequalities for  Young inequality was given by M. Tominaga in \cite{Tom}, 
using the Specht ratio, in the following way
\begin{equation}\label{eq_1_2}
 (1-\lambda) a+ \lambda b\leq   S\left(\frac{a}{b}\right)a^{1-\lambda} b^{\lambda}, 
\end{equation}
for positive real numbers $a,b$ and $\lambda \in[0,1]$, where the
Specht ratio \cite{Spe,Furuta,FMPS}, was defined by
$$S(h)\equiv \frac{h^{\frac{1}{h-1}}}{e\log h^{\frac{1}{h-1}}}, \ (h\not=1)$$
for a positive real number $h$. Note that $\lim_{h\to 1} S(h) =1$ and $S(h)=S(1/h) > 1$ for $h\neq 1, h>0$.
We call the inequality (\ref{eq_1_2}) a ratio-type reverse inequalitiy for Young's inequality. 
M. Tominaga also showed in \cite{Tom} the following inequality:
\begin{equation}\label{eq_1_3}
(1-\lambda) a+ \lambda b\leq  L(a,b) \log S\left(\frac{a}{b}\right) +a^{1-\lambda} b^{\lambda}, 
\end{equation}
for positive real numbers $a,b$ and $\lambda \in[0,1]$, where the logarithmic mean \cite{Kubo} $L(x,y)$ is defined by
$$
L(x,y) \equiv \frac{y-x}{\log y-\log x},\,\,(x \neq y),\quad L(x,x)=x.
$$
 We call the inequality (\ref{eq_1_3}) a difference-type reverse inequalitiy for Young's inequality. 
Based on the scalar inequalities (\ref{eq_1_2}) and  (\ref{eq_1_3}), M. Tominaga showed the following two reverse inequalities
for invertible positive operators.

\begin{The}{\bf (\cite{Tom})}
For invertible positive operators $A$ and $B$ with $0<mI\leq A,B\leq MI$, we have
\begin{itemize}
\item[(i)] (Ratio-type reverse inequality)
\begin{equation} \label{eq_1_4}
(1-\lambda)A+\lambda B \leq S(h) A\sharp_{\lambda} B. 
\end{equation}
\item[(ii)] (Difference-type reverse inequality)
\begin{equation} \label{eq_1_5}
(1-\lambda)A+\lambda B \leq  A\sharp_{\lambda} B +L(1,h) \log S(h) B. 
\end{equation}
\end{itemize}
\end{The}

Our purpose of this short paper is to give two reverse inequalities which are different from (\ref{eq_1_4}) and  (\ref{eq_1_5}).
 

\section{Main results}
We first show the following remarkable scalar inequality.
\begin{The} \label{the_2_1}
Let $f:[a,b]\to\mathbb{R}$ be a twice differentiable function such that
there exist real constant $M$ so that $0\leq f^{\prime\prime}(x)\leq
M$ for $x\in[a,b]$.\\
Then we have
\begin{equation} \label{eq_2_1}
0\leq (1-\lambda) f(a)+ \lambda f(b)-f((1-\lambda) a+ \lambda b)\leq
\lambda(1-\lambda)M(b-a)^2,
\end{equation}
where $\lambda\in[0,1]$.
\end{The}

{\it Proof:} The first part of inequality (\ref{eq_2_1}) holds because $f$ is a convex function. 
Next, we prove second part of inequality (\ref{eq_2_1}). 

For $\lambda\in\{0,1\}$, we obtain the equality in relation  (\ref{eq_2_1}). 
Now, we consider $\lambda\in(0,1)$, which means that $a<(1-\lambda)
a+ \lambda b<b$ and we use Lagrange's theorem for the intervals
$[a,(1-\lambda) a+ \lambda b]$ and $[(1-\lambda) a+ \lambda b,b]$.
Therefore, there exists real numbers $c_1\in(a,(1-\lambda)
a+ \lambda b)$ and $c_2\in((1-\lambda) a+ \lambda b,b)$ such that
\begin{equation} \label{eq_2_2}
f((1-\lambda) a+ \lambda b)-f(a)= \lambda (b-a)f^\prime(c_1)
\end{equation}
and
\begin{equation}\label{eq_2_3}
f(b)-f((1-\lambda) a+ \lambda b)=(1-\lambda)(b-a)f^\prime(c_2).
\end{equation}
Multiplying relation (\ref{eq_2_2}) by $(1-\lambda)$ and relation $(\ref{eq_2_3})$ by
$ \lambda$, and then adding, we deduce the following relation:
$$
(1-\lambda) f(a)+ \lambda f(b)-f((1-\lambda) a+ \lambda b) 
=\lambda (1-\lambda)(b-a)[f^\prime (c_2)-f^\prime(c_1)].
$$
Again, applying Lagrange's theorem on the interval $[c_1,c_2]$, we
obtain
\begin{equation}\label{eq_2_4}
(1-\lambda) f(a)+ \lambda f(b)-f((1-\lambda) a+ \lambda b) 
=\lambda(1-\lambda)(b-a)(c_2-c_1)f^{\prime\prime}(c), 
\end{equation}
where $c\in(c_1,c_2)$.
Since $0\leq f^{\prime\prime}(x)\leq M$ for $x\in[a,b]$ and $c_2-c_1\leq b-a$ and
making use of relation (\ref{eq_2_4}), we obtain the inequality (\ref{eq_2_1}).

\hfill \qed

\begin{Cor}
For $a, b>0$ and $\lambda\in[0,1]$, the following inequalities hold.
\begin{itemize}
\item[(i)] (Ratio-type reverse inequality)
\begin{equation} \label{eq_2_5}
a^{1-\lambda} b^{\lambda}\leq (1-\lambda) a+ \lambda b\leq a^{1-\lambda}
b^{\lambda}\exp\left\{ \frac{\lambda(1-\lambda)(a-b)^2}{d_1^2} \right\},
\end{equation}
where $d_1\equiv \min\left\{a,b\right\}$.
\item[(ii)] (Difference-type reverse inequality)
\begin{equation}\label{eq_2_6}
a^{1-\lambda} b^{\lambda}   \leq (1-\lambda) a+ \lambda b  \leq a^{1-\lambda} b^{\lambda}   +\lambda
(1-\lambda)\left\{ \log\left( \frac{a}{b} \right) \right\}^2 d_2,
\end{equation}
where  $d_2    \equiv \max\left\{a,b \right\}$.
\end{itemize}
\end{Cor}

{\it Proof:} 
\begin{itemize}
\item[(i)]
It is easy to see that if we take $f(x)=-\log x$ in Theorem \ref{the_2_1}, then we have
$$
\log\left\{(1-\lambda)a +\lambda b\right\} \leq \log \left(a^{1-\lambda} b^{\lambda}\right) +\log\left( \exp\left\{\lambda(1-\lambda)  
f^{\prime\prime}(c)  (b-a)^2\right\}\right)
$$
which implies  inequality (\ref{eq_2_5}), since $f^{\prime\prime}(c) = \frac{1}{c^2} \geq \frac{1}{d_1^2}$. 
\item[(ii)]
If we take $f(x)=e^x$ (which is convex on $(-\infty,\infty)$) in Theorem \ref{the_2_1}, we
 obtain the relation
$$
0\leq \lambda e^{\alpha}+(1-\lambda)e^{\beta}-e^{\lambda \alpha+(1-\lambda)\beta}\leq
\lambda(1-\lambda) ( \alpha-\beta)^2  f^{\prime\prime}(\gamma),
$$
where $\gamma \equiv \max\left\{\alpha,\beta\right\}$.
Putting $a = e^{\alpha}$ and $b=e^{\beta}$, then we have
$$
0\leq (1-\lambda)a+\lambda b- a^{1-\lambda}b^{\lambda} \leq \lambda(1-\lambda)e^c \left(\log \frac{a}{b}\right)^2
$$
where $c\equiv \max \left\{\log a,\log b\right\}$.
Thus we have inequality (\ref{eq_2_6}), putting $d_2=e^c$.

\end{itemize}

\hfill \qed

From here, we consider bounded linear operators acting on a complex Hilbert space $\fH$. 
If a bounded linear operator $A$ satisfies $A=A^*$, then $A$ is called a self-adjoint operator. 
If a self-adjoint operator $A$ satisfies $\langle x \vert A \vert x \rangle \geq 0$ for any $\vert x \rangle \in \fH$,
then $A$ is called a positive operator. In addition, $A \geq B$ means $A -B \geq 0$.

\begin{The}
For $\lambda \in [0,1]$, two invertible positive operators $A$ and $B$ satisfying the ordering $mI \leq A \leq B\leq M I \leq I$ with $h\equiv \frac{M}{m}$, 
we have the following operaor inequalities.
\begin{itemize}
\item[(i)] (Ratio-type reverse inequality)
\begin{equation}\label{eq_the_3_01}
A\sharp_{\lambda}B \leq (1-\lambda) A+\lambda B \leq \exp\left(\lambda(1-\lambda)\left(1-h\right)^2\right) A\sharp_{\lambda}B. 
\end{equation}
\item[(ii)] (Difference-type reverse inequality)
\begin{equation}   \label{eq_the_3_02}
A\sharp_{\lambda} B \leq (1-\lambda) A+\lambda B  \leq A\sharp_{\lambda} B +\lambda (1-\lambda) \left( \log h \right)^2 B. 
\end{equation}
\end{itemize}
\end{The}

{\it Proof:} 
\begin{itemize}
\item[(i)] 
The first inequalities in (\ref{eq_the_3_01}) and  (\ref{eq_the_3_02}) are well-known so that 
we prove the two second inequalities in (\ref{eq_the_3_01}) and  (\ref{eq_the_3_02}). 
From the inequality (\ref{eq_2_5}) with $a \leq b$, we have 
$$
(1-\lambda)t +\lambda \leq t^{1-\lambda}  e^{\lambda(1-\lambda)(1-\frac{1}{t})^2},
$$
for $0<t \leq 1$. Thus we have the following inequality for the invertible positive operator $mI \leq T\leq MI\leq I$:
$$
(1-\lambda) T+\lambda \leq T^{1-\lambda} \max_{m\leq t\leq M}  e^{\lambda(1-\lambda)(1-\frac{1}{t})^2}.
$$
Putting $T\equiv B^{-1/2}AB^{-1/2} \leq I$ (which satisfies $A\leq B$), then we have $\frac{1}{h}\leq B^{-1/2}AB^{-1/2} \leq h$, and then
 we have
$$
(1-\lambda) B^{-1/2}AB^{-1/2} +\lambda \leq \left( B^{-1/2}AB^{-1/2} \right)^{1-\lambda}  \max_{\frac{1}{h}\leq t\le 1 \leq h}    e^{\lambda(1-\lambda)(1-\frac{1}{t})^2}.
$$
Multiplying $B^{1/2}$ to the both sides in the above inequality, we obtain the inequality (\ref{eq_the_3_01}), since $A\sharp_{\lambda}B =B\sharp_{1-\lambda}A$.

\item[(ii)] 
By the similar way to the proof of the second inequality in (\ref{eq_the_3_01}) from the inequality (\ref{eq_2_6}) with $1 \leq a \leq b$,  we have 
$$
(1-\lambda) t +\lambda -t^{1-\lambda} \leq \lambda (1-\lambda)   \left( \log t \right)^2,
$$
for $0<t \leq 1$. Thus we have the following inequality for the invertible positive operator  $mI \leq T\leq MI\leq I$:
$$
(1-\lambda) T +\lambda -T^{1-\lambda} \leq \lambda (1-\lambda)  \max_{m\leq t\leq M}   \left( \log t \right)^2.
$$
Putting $T\equiv B^{-1/2}AB^{-1/2} \leq I$ (which satisfies $A\leq B$), then we have $\frac{1}{h}\leq B^{-1/2}AB^{-1/2} \leq h$, and then we have
$$
(1-\lambda) B^{-1/2}AB^{-1/2} +\lambda -\left(B^{-1/2}AB^{-1/2} \right)^{1-\lambda} \leq \lambda (1-\lambda) \max_{\frac{1}{h}\leq t\leq h}    \left( \log t \right)^2,
$$
which implies the inequality (\ref{eq_the_3_02}), by multiplying $B^{1/2}$ to the both sides in the above inequality,  since $A\sharp_{\lambda}B =B\sharp_{1-\lambda}A$.
\end{itemize}

\hfill \qed

\begin{Rem}
It is natural to consider that our inequalities are better than Tominaga's inequalities under the assumpution $A\leq B$.
Firstly we compare our inequality (\ref{eq_2_5}) with (\ref{eq_1_2}). 
For this purpose we take two numerical example under the condition $0<t \leq 1$.
\begin{itemize}
\item[(i)] Take $t=\frac{1}{2}$ and $\lambda =\frac{1}{20}$, then we have
$$
\exp\left(\lambda(1-\lambda)\left(1-\frac{1}{t}\right)^2 \right)-S(t) \simeq -0.0128295.
$$
\item[(ii)] Take $t=\frac{1}{2}$ and $\lambda =\frac{1}{10}$, then we have
$$
\exp\left(\lambda(1-\lambda)\left(1-\frac{1}{t}\right)^2 \right)-S(t) \simeq 0.0326986.
$$
\end{itemize}
Thus we can conclude that there is no ordering between  (\ref{eq_2_5}) and (\ref{eq_1_2}).
If we compare our factor in the right hand side of the inequality (\ref{eq_the_3_01}) with one in (\ref{eq_1_4}), then the similar result follows, putting $h=2 >1$ with  $\lambda =\frac{1}{20}$ or $\lambda =\frac{1}{10}$, for example.

Secondly we compare our inequality  (\ref{eq_2_6}) with (\ref{eq_1_3}). 
For this purpose we take two numerical example under the condition $0<t \leq 1$.
\begin{itemize}
\item[(i)] Take $t=\frac{1}{2}$ and $\lambda =\frac{1}{5}$, then we have
$$
L(1,t)\log S(t)-\lambda(1-\lambda)\left(\log t\right)^2 \simeq -0.0338368.
$$
\item[(ii)] Take $t=\frac{1}{2}$ and $\lambda =\frac{1}{20}$, then we have
$$
L(1,t)\log S(t)-\lambda(1-\lambda)\left(\log t\right)^2 \simeq 0.0202141.
$$
\end{itemize}
Thus we can conclude that there is no ordering between  (\ref{eq_2_6}) and (\ref{eq_1_3}).

Therefore we may conclude our two reverse inequalities for Young's inequality do not trivial results under the assumpution $A \leq B$. 

\end{Rem}

\section*{Acknowledgement}
The author (S.F.) was supported in part by the Japanese Ministry of Education, Science, Sports and Culture, Grant-in-Aid for
Encouragement of Young Scientists (B), 20740067
The author (N.M.) was supported in part by the Romanian Ministry of Education, Research
and Innovation through the PNII Idei project 842/2008.

\end{document}